\documentclass[a4paper,12pt]{article}
\usepackage[T1]{fontenc}
\usepackage[cp1250]{inputenc}
\usepackage[pdftex]{graphicx}
\usepackage{amsmath,amstext,amssymb,amsopn,amsthm,amsfonts,mathrsfs,extraipa,dsfont,bbm}

\newcommand{\A}{\mathcal{A}}
\newcommand{\B}{\mathcal{B}}
\newcommand{\I}{\mathcal{I}}

\newcommand{\ccc}{\mathfrak{c}}
\newcommand{\F}{\mathfrak{F}}
\newcommand{\R}{\mathbb{R}}

\newcommand{\nsk}{[\omega]^\omega}
\newcommand{\sk}{[\omega]^{<\omega}}

\newcommand{\el}{\nsk_{\texttt{\tiny EL}}}

\begin{document}

\title{{\bf On partitions of Ellentuck-large sets}}
\author{Ryszard Frankiewicz \and S{\l}awomir Szczepaniak}
\date{}
\maketitle

\begin{abstract}
It is proved that no non-meager subspace of the space $\nsk$ equipped
with the Ellentuck topology does admit a Kuratowski partition,
that is such a subset cannot be covered by a family $\F$ of disjoint relatively meager sets
such that $\bigcup\F'$ has the Baire property (relatively) for every
subfamily $\F'\subseteq\F$. Some remarks concerning continuous restrictions of functions with domain in the Ellentuck space are made.
\end{abstract}

\vspace{1cm}

We consider the set $\nsk$ of infinite subsets of $\omega$ equipped with Ellentuck topology. The sets of the form
$[a,A]=\{B\in[A]^\omega: a\sqsubset B\subseteq a\cup A\}$, where
$a\in\sk$ and $A\in\nsk$, establish a base of this topology. Through the whole paper we called them {\em basis sets} (the
expression ''$a\sqsubset B$'' stands for {\em $a$ is a initial segment of
$B$}). For unexplained topological notions we refer to \cite{En} or \cite{Ku}.
This space, called here {\em the Ellentuck space} and denoted by $\el$, is quite well studied. For
example it is widely known that $\el$ is neither compact nor
metrizable (it is even non-normal, see \cite{P}), but it is a Baire space.
Moreover, its nowhere dense sets form a $\sigma$-ideal \cite{T}.
These nowhere dense sets $X\subseteq\el$, called {\em Ramsey null sets}, are characterized as follows: for every set of form $[a,A]$ there exists $B\in[A]^\omega$ such that $[a,B]\cap X=\emptyset$. The $\sigma$-ideal of Ramsey null is a subfamily of $\sigma$-algebra
of so called {\em completely Ramsey} sets. We say that $X\subseteq\el$ is completely Ramsey if for every set of form $[a,A]$ there exists $B\in[A]^\omega$ such that either $[a,B]\subseteq X$ or $[a,B]\cap X=\emptyset$. It turns out that completely Ramsey sets in $\el$ are exactly the sets having the Baire property in the Ellentuck topology (\cite{T}). The above nice combinatorial characterizations are what make $\el$ close to more familiar Polish spaces. For example, the following was proved in \cite{LS} (cf. \cite{AT})

\medskip

\noindent----------------------

\noindent{\small 2010 {\em Mathematics Subject Classification}. 54B15 (Primary), 54D20 (Secondary).

}

\medskip

\noindent{\bf  Louveau-Simpson Theorem} {\em Let $\F$ be a point finite family (i.e. the same element can appear in members of $\F$ only finitely many times) of Ramsey null sets such that the union of any subfamily is completely Ramsey. Then the union of the whole family is Ramsey null.}

\medskip

\noindent The above theorem as well as its proof strikingly resembles analogous results of Solovay, Prikry and Bukovsk\'{y} (cf. \cite{FG}) which were generalized to

%
%
%
%

\medskip

\noindent{\bf Four Poles Theorem \cite{BCGR}} {\em Let $\B$ be the $\sigma$-algebra generated by Borel sets in the Polish space $X$ and a $\sigma$-ideal $\I$ on $X$ with the Borel basis. Then any point finite family of sets from $\I$ which covers $X$ has a subfamily with not $\B$-measurable union.}

\medskip

\noindent Actually, as noted in \cite{FG}, Louveau-Simpson Theorem was known earlier as an instance of some more general results from \cite{EFK2} concerning so called {\em pseudobasically compact} spaces since $\el$ is one of them (\cite{FG}). From this point of view one can see that Ellentuck space also shares common features with compact spaces.

The presented paper deals with a strictly topological version of conclusion of Four Poles Theorem. Therefore let us make the following

\medskip

\noindent{\bf Definition} Let $(X,\tau)$ be a topological space and let $\F$ be a partition of $X$ into meager sets.
We say that $\F$ is a {\em Kuratowski partition} if $\bigcup\F'$ has the Baire property for any subfamily $\F'$ of $\F$.

\medskip

\noindent By replacing in the above definition "partition" by "point finite family" we usually gain no new results nor we lose old ones. In this terminology Louveau-Simpson Theorem states that no non-meager subspace of $\el$ with the Baire property admits a Kuratowski partition. Let us now argument that a Kuratowski partition is not only a technical notion, its main motivation is hidden in the following equivalence

\medskip

\noindent{\bf Proposition 1 (\cite{EFK1}, \cite{FK})} {\em For a topological space $X$ and $Y$, where the last one possesses $\sigma$-disjoint base, the following are equivalent:

$(i)$ No subspace of $X$ of form $G\setminus F$, with $G$ open and $F$ meager, does admit a Kuratowski partition.

$(ii)$ Every Baire-measurable function $f:X\mapsto Y$ is continuous on a co-meager subset of $X$.

The same equivalence is obtained if $Y$ is replaced by a metric space or even the one with the discrete topology.}

\medskip

\noindent The above means that Kuratowski partitions serves as tools in studying the following problem: {\em to which extent one can realize in topological context a widely known Luzin's theorem on measurable functions?} To present briefly a story of this problem let us make the following ad hoc definition. Call a pair of topological spaces $(X,Y)$ {\em a Luzin pair} if $(ii)$ from {\bf Proposition 1} is satisfied. Firstly O.Nikodym in 1929 established that Luzin pairs exists; he proved actually that $(\R,\R)$ is a Luzin pair (\cite{N}). Several months later K.Kuratowski (\cite{Ko}) noted that Nikodym's proof works also for any pairs $(X,Y)$ with second countable $Y$ and arbitrary $X$. In 1935 (\cite{Kp}) K.Kuratowski asked if any pair $(X,Y)$ with $X$ being completely metrizable and $Y$ being arbitrary metric space is a Luzin pair. This was answered affirmatively in \cite{EFK1} for $X$ with weight $\leqslant\ccc$ and negatively in general in \cite{FK}, where the following was proved.

\medskip

\noindent{\bf Proposition 2} {\em The following theories are equiconsistent

$(i)$ ZFC+"there exists a measurable cardinal"

$(ii)$ ZFC+"some complete metric space $X$ admits Kuratowski partition"

$(iii)$ ZFC+"some Baire metric space $X$ admits Kuratowski partition"}

\medskip

\noindent Therefore the quest for Luzin pairs $(X,Y)$ makes sense only when $Y$ is like in {\bf Proposition 1}. By the same proposition it is reduced to the question about Kuratowski partitions of large (in sense of category, i.e. non-meager) subspace of a (Baire) metric space $X$. We shall show that the aforementioned problem can be reduced even more.

\medskip

\noindent{\bf Proposition 3} {\em If there exists a Baire space with a Kuratowski partition,
then exists a Baire metric space admitting a Kuratowski partition. In particular the existence of Kuratowski partition is not a metric problem.}

\medskip

\noindent We postpone the proof of {\bf Proposition 3} to the end of the paper after clarifying its purpose in more details.

The space $X$ from {\bf Proposition 2} was constructed as a subspace of a space $^\omega(2^{\omega_1})$ with the standard product topology. Much earlier it was also known that whole $^\omega(2^{\omega_1})$ does not admit a Kuratowski partition, thus sharing this property with $\el$.
Therefore, some time ago the question arose whether there exists (at least consistently) a subspace of $\el$ with a Kuratowski partition and whether consistency strength of its existence is less than the one in {\bf Proposition 2}. It was hoped (by Shelah among many people) that nice combinatorial properties of $\el$ and its resemblance to metrizable and compact spaces should provide such a 'natural' example of a space with a Kuratowski partition. Furthermore, it follows from {\bf Proposition 3} that aforementioned 'naturalness' might be hoped to decrease consistency strength from {\bf Proposition 2 $(i)$} and could serve as a standard example of such object. We now show that the above hopes
were vain and the hypothesis was false.


\newpage

\noindent{\bf Proposition 4}\,{\em\,No\,non-meager\,subspace\,of\,$\el$\,admits\,a\,Kuratowski\,partition.}

\medskip

\noindent{\bf Proof:} Let $X$ be a non-meager subspace of $\el$ and denote by $\B$ a family of basis sets in $\el$. Not being Ramsey null (RN in short) $X$ is dense in some $B\in\B$ and since $B$ is homeomorphic to $\el$ we may assume w.l.o.g. that $X$ is already dense in $\el$. Let $\F$ be a partition of $X$ into RN-sets. To prove the theorem we shall find $\F'\subseteq\F$ such that $\bigcup\F'$ lacks the Baire property in $X$. Since $|\F|\leqslant |2^\omega|$ enumerate $\F$ by $E\subseteq2^\omega$, putting $F_f:=\emptyset$ for $f\in 2^\omega\setminus E$. Next, for all $s\in 2^{<\omega}$ define $F[s]:=\bigcup\left\{F_f\in\F:s\subseteq f\right\}$. Note that $\left\{F[s]:s\in 2^{<\omega}\right\}$ forms a Cantor scheme on $X$ (cf. \cite{K}). Moreover, a tree $\mathcal{S}:=\left\{s\in 2^{<\omega}:F[s]\,\,\mbox{is not RN}\right\}$ is a perfect subtree of $2^{<\omega}$. Indeed, for if it is not the case one can find $s\in\mathcal{S}$ such that for each $t\supseteq s$ there is $i<2$ with $F[t\,\hat{}\,i]$ being a RN-set. Then a tree $\left\{F[s\,\hat{}\,t]:t\in 2^{<\omega}\right\}$ has only one branch indexed by $f_s\in E$ consisting of non-meager sets and for all $f\notin E$ there is a RN-set $F[s_f]$, $s_f\subseteq f$, being a superset of $F_f$. Thus, a non-RN-set $F[s]$ can be covered by countably many RN-sets ($F_{f_s}$ and $F[s]$ for $s\nsubseteq f_s$), a contradiction. Similarly one proves that $X':=\bigcup\left\{F_f\in\F:f\in[\mathcal{S}]\right\}$ is a comeager subset of $X$. Hence, $\F$ is a Kuratowski partition of $X$ iff
$\left\{F_f\in\F:f\in[\mathcal{S}]\right\}$ is a Kuratowski partition of $X'$. So assume that $F[s]$ is not RN for all $s\in 2^{<\omega}$.

Define $G_s:=\bigcup\left\{B\in\B:F[s]\,\,\mbox{is dense in}\,\,B\right\}$ for $s\in 2^{<\omega}$. We claim there is $s\in 2^{<\omega}$ with $G_{s\,\hat{}\,0}\cap G_{s\,\hat{}\,1}\neq\emptyset$. Then any $B\in\B$ included in this intersection witnesses that $F[s\,\hat{}\,0]$ lacks the Baire property in $X$. Towards a contradiction assume otherwise. Therefore a family $\mathcal{U}:=\{G_s:s\in 2^{<\omega}\}$ forms a Cantor scheme of open sets on $\el$. Put $\mathcal{G}(C):=\left\{\bigcap_nG_{f\upharpoonright n}:f\in C\right\}\setminus\{\emptyset\}$ for any $C\subseteq2^\omega$. Note that a set $G(n):=\bigcup\{G_s:s\in 2^n\}$, $n<\omega$, is dense in $X$. Thus $\bigcup\mathcal{G}(2^\omega)=\bigcap_nG(n)$ is a dense $G_\delta$-set in $\el$ so it has a dense interior by Th.3.1. of \cite{Hr}. Hence $F:=\el\setminus\bigcup\mathcal{G}(2^\omega)$ is a RN-set. We now show

\medskip

\noindent{\bf Claim} {\em A family $\displaystyle\mathcal{G}:=\mathcal{G}(2^\omega)\cup\left\{F\right\}\setminus\{\emptyset\}$ is a Kuratowski partition of $\el$.}

\medskip

Indeed, the family $\mathcal{G}$ is a partition by the fact that $\mathcal{U}$ forms a Cantor scheme on $\el$. Moreover, any element of $\mathcal{G}$ is RN. This was already shown for $F$ so consider $f\in E$. As $\bigcap_nG_{f\upharpoonright n}$ is a $G_\delta$-set in $\el$, by Th.3.1. of \cite{Hr}, it suffices to prove that it has an empty interior in $\el$. If however there is $B\in\B$ with $B\subseteq\bigcap_nG_{f\upharpoonright n}$, then, by the definition of $G_s$, it would hold $B\subseteq\bigcap_n \overline{F\left[f\upharpoonright n\right]}=\overline{\bigcap_n F\left[f\upharpoonright n\right]}=\overline{F_f}$ implying that $F_f$ is not a RN-set, a contradiction. It remains to show that for any $\mathcal{G}'\subseteq\mathcal{G}$ a set $\bigcup\mathcal{G}'$ has the Baire property in $\el$. We need only to consider a subfamily of $\mathcal{G}\setminus\left\{F\right\}$. As any such family is of the form $\mathcal{G}(C)$ for some $C\subseteq E$ we obtain $\bigcup\mathcal{G}(C)=\bigcup_{f\in C}\bigcap_nG_{f\upharpoonright n}=
\bigcap_n\left\{G_{f\upharpoonright n}:f\in C\right\}$. The second equality holds as $\mathcal{U}$ is a Cantor scheme. Therefore $\bigcup\mathcal{G}(C)$, being $G_\delta$-set, has the Baire property in $\el$. This ends the proof of {\bf Proposition 4} as Claim contradicts Louveau-Simpson Theorem.\quad$\Box$

\medskip

\noindent{\bf Remark} The problem of an existence of a 'small' Baire spaces (preferably associated with some classical forcing notion) possessing a Kuratowski partition remains open.

\medskip

\noindent{\bf Remark} Besides Luzin-Nikodym type theorem there is one more theorem in Analysis concerning continuous restrictions of functions. This is Blumberg Theorem from 1922 \cite{B}. Blumberg proved that any function $f:\R\to\R$ is continuous on dense subset of $\R$. This was further generalized by J.C.Bradford and C.Goffman in \cite{BG} to real-valued functions defined on an arbitrary Baire metric space. In fact their proof works also for arbitrary Baire spaces as noted in \cite{W}. Hence any function $f:\el\to\R$ is continuous on dense subset of $\el$. In fact, by combining results from Proposition 1.4 \cite{W} and Theorem 1 \cite{P}, space $\R$ can be replaced by any topological space with weight less than distributivity number of the Boolean algebra $P(\omega)/Fin$ (cf. \cite{BPS}) which is always between $\omega_1$ and $\ccc$.

There exist strengthenings and limitations for classical Blumberg Theorem. We say that subset of topological space $X$ is categorically dense or nowhere meager (resp. $\kappa$-dense) if it meets any open subset in a non-meager set (resp. in a set of cardinality $\kappa$) \cite{Br}. These two kinds of densities can be viewed as the strongest in topological (resp. set-theoretical) sense. In the Ellentuck space however any dense sets are dense in these two senses which trivializes the problem of strengthening Blumberg Theorem for the Ellentuck space. Indeed, since meager subsets of $\el$ are nowhere dense, any dense subset is categorically dense. Moreover, dense sets are $\ccc$-dense as any non-meager subset of $\el$ is of size $\ccc$. For if $X$ is non-meager it is dense in some basis set $[a,A]$. Let $\A$ be an almost disjoint family of size $\ccc$ of subsets of $A\setminus\max a$. Then $X$ meets any member of pairwise disjoint family $\{[a,B]:B\in\A\}$.

As far as limitations are concerned remind that the well-known obstruction to Blumberg Theorem is a Sierpinski-Zygmund function $f:\R\to\R$ which is discontinuous on any set of reals of size $\ccc$. As mentioned above the Ellentuck space satisfies strong form of Blumberg Theorem thus any reasonable candidate for Sierpinski-Zygmund-like function turns out to be trivial. For example a discontinuous functions on any basis set could be such a candidate. However it easy to find such an example: just take a characteristic function of a Bernstein set in standard metrizable topology of $\nsk$ (\cite{T}, cf. \cite{P}). Such a set splits any basis set into two disjoint ones therefore its characteristic function does not even possess Baire-measurable restriction to any basis set.

\medskip

We finish the paper with the promised proof of {\em Proposition 3}. First we need some preparations. For topological space $(X,\tau)$ denote $\tau^+=\tau\setminus\{\emptyset\}$. Define
$$X(\tau):=\left\{x\in(\tau^+)^{\omega}:\,\bigcap_{n<\omega}x(n)\neq\emptyset\right\}.$$
Treated as subspace of a complete metric space
$^\omega(\tau^+)$, where the set $\tau^+$ is equipped with the
discrete topology, the space $X(\tau)$ is a metric space. A basis
of the space $X(\tau)$ is given by the sets of the form
$[s]:=\{x\in X(\tau):x\supseteq s\}$ for $s\in(\tau^+)^{<\omega}$. After the paper was finished we learned from Piotr Zakrzewski that a version of the space $X(\tau)$ was defined earlier in \cite{KR}, where the author also proved generalization of {\bf Lemma 4}. In this lemma we shall prove that $X(\tau)$ is a Baire space provided that $(X,\tau)$ is a Baire space. We use the following well-known theorem due to Oxtoby

\medskip

\noindent{\bf Oxtoby Theorem (\cite{K}, 8.11)} {\em A nonempty topological space $(X,\tau)$ is a Baire space iff player I has no winning strategy in the Choquet game $\mathcal{G}(X)$.}

\medskip

\noindent Recall only that Choquet game $\mathcal{G}(X)$ of $X$ consists in alternating choices (made by players I and II) of nonempty open sets in $X$. I player starts and play with open sets $U_n$'s and II player responds with open sets $V_n$'s in such a way that $U_0\supseteq V_0\supseteq U_1\supseteq V_1\supseteq\dots$. Player I (II) wins a run $(U_0,V_0,U_1,V_1,\dots)$ of the game $\mathcal{G}(X)$ if $\bigcap_nU_n=\bigcap_nV_n=\emptyset$ ($\neq\emptyset$). For other undefined notions concerning topological games (a winning strategy, an equivalent game, etc.) we refer to \cite{K} (8.10, 8.36). We only remark that if in the above definition players are allowed only to choose open sets from a fixed basis then this modified game is equivalent to the Choquet game $\mathcal{G}(X)$ of $X$. In the case of $X(\tau)$ this means actually that both players made their moves in the modified Choquet game of $X(\tau)$ by extending finite sequences of elements of $\tau^+$ chosen by the second player.

\medskip

\noindent{\bf Lemma 4} {\em If $(X,\tau)$ is a Baire space then so is $X(\tau)$.}

\medskip

\noindent{\bf Proof:}

\medskip

Towards a contradiction suppose $X(\tau)$ is not a Baire space. Therefore by Oxtoby's theorem player I has a winning strategy in the (modified) Choquet game $\mathcal{G}(X(\tau))$. We shall describe a winning strategy for player I in the Choquet game $\mathcal{G}(X)$ which contradicts (again via Oxtoby's theorem) that $(X,\tau)$ is a Baire space. Put $U_0=\bigcap\{s_0(k):k\in\mbox{dom}(s_0)\}$ where $s_0\in(\tau^+)^{<\omega}$ is a first move of player I in some fixed winning strategy for I in the game $\mathcal{G}(X(\tau))$. Note that $U_0$ is nonempty. Indeed since $[s_0]\neq\emptyset$ (by the definition of moves in Choquet games) there is $x_0\in X(\tau)$ such that $x_0\supseteq s_0$; so $U_0\supseteq\bigcap_nx_0(n)\neq\emptyset$. Let $n<\omega$ and suppose that player II responds with $V_n\subseteq U_n$ to the $n+1^{\mbox{\tiny th}}$ move of player I. Define $s_{n+1}\in(\tau^+)^{<\omega}$ as the unique extension (response to II's move) of a sequence $s_n\,\hat{}\,(V_n)$ in the I's winning strategy in the game $\mathcal{G}(X(\tau))$. Put $U_{n+1}=\bigcap\{s_{n+1}(k):k\in\mbox{dom}(s_{n+1})\}$ and notice it is a nonempty set by the same reason as $U_0$ was.

Put $x:=\bigcup_ns_n\in(\tau^+)^{\omega}$. Observe that $\bigcap_n[s_n]=\emptyset$ since a sequence $(s_n)_n$ is a run of player I in his winning strategy in the game $\mathcal{G}(X(\tau))$. This means that $\bigcap_{k<\omega}x(k)=\emptyset$ by the definitions of $X(\tau)$ and its basic open sets. Therefore $$\bigcap_{n<\omega}U_n=\bigcap_{n<\omega}\bigcap_{k\in\mbox{\tiny dom}(s_n)}s_n(k)=\bigcap_{k<\omega}x(k)=\emptyset.$$
Hence the strategy for player I described above is the winning one. \quad $\Box$

\medskip


\noindent{\bf Lemma 5} {\em If $(X,\tau)$ does admit a Kuratowski partition, then so does $X(\tau)$.}

\medskip

\noindent{\bf Proof:}

\medskip

Fix a Kuratowski partition $\F$ of $X$. Define a function $\varphi:X(\tau)\to X$ by $$\varphi(x)=\min\bigcap_{n<\omega}x(n),\quad x\in X(\tau),$$
where $\min$ refers to a minimum with respect to some fixed well-ordering of $X$. In fact any function with $\varphi(x)\in\bigcap_nx(n)$ for $x\in X(\tau)$ could serve for our purposes. Let us verify some properties of $\varphi$.

\begin{enumerate}
  \item {\em If $N$ is nowhere dense in $X$ then $\varphi^{-1}[N]$ is nowhere dense in $X(\tau)$.}

  Indeed, let $s\in(\tau^+)^{<\omega}$ be arbitrary such that $[s]\neq\emptyset$; in particular $U:=\bigcap\{s(k):k\in\mbox{dom}(s)\}\neq\emptyset$ (see the argument for nonemptyness of $U_0$ in {\bf Lemma 3.2}). Since $N$ is nowhere dense in $X$ let $V\in\tau^+$ be such that $V\subseteq U$ and $V\cap N=\emptyset$. Put $t:=s\,\hat{}\,(V)$ and observe $[t]\neq\emptyset$ as $\bigcap\{t(k):k\in\mbox{dom}(t)\}=U\cap V=V$. We need only to check that $[t]\cap\varphi^{-1}[N]=\emptyset$. Toward a contradiction suppose that there is $x\in [t]$ with $\varphi(x)\in N$. However from $x\supseteq t$ it follows that $\varphi(x)\in\bigcap_nx(n)\subseteq\bigcap\{t(k):k\in\mbox{dom}(t)\}=V$. Hence $\varphi(x)\in V\cap N$ contradicting the choice of $V$.
  \item {\em If $M$ is meager in $X$ then $\varphi^{-1}[M]$ is meager in $X(\tau)$.}

  Indeed, let $(N_n)_n$ be a sequence of nowhere dense sets in $X$ such that $M=\bigcup_nN_n$. Then $\varphi^{-1}[M]=\bigcup_n\varphi^{-1}[N_n]$ and by $1.$ for all $n<\omega$ the set $\varphi^{-1}[N_n]$ is nowhere dense in $X(\tau)$.
  \item {\em The function $\varphi:X(\tau)\to X$ is Baire-measurable.}


  Indeed, let $U\in\tau^+$. Then the following set is open in $X(\tau)$ $$\mathcal{U}:=\left\{x\in X(\tau):\mbox{rng}(x)\cap P(U)\cap\tau^+\neq\emptyset\right\}=\bigcup_{\tiny U\supseteq V\in\tau^+}\bigcup_{n<\omega}\left\{x\in X(\tau):x(n)=V\right\}.$$ Now $\varphi(x)\in\bigcap_nx(n)$ for $x\in X(\tau)$ and $\bigcap_nx(n)\subseteq U$ if $U\in\mbox{rng}(x)$ imply $$\varphi[\mathcal{U}]\subseteq\bigcup\left\{\bigcap_nx(n)\in X:U\in\mbox{rng}(x)\right\}\subseteq U.$$ Hence it suffices to show that $\mathcal{U}$ is dense in $\varphi^{-1}[U]$. Towards a contradiction suppose that $x\in\varphi^{-1}[U]\setminus\mbox{cl}\,\mathcal{U}$ for some $x\in X(\tau)$. Let $[s]$ be a basic open neighbourhood of $x$ in $X(\tau)$ omitting $\mathcal{U}$. Observe that $V:=\bigcap\{s(k):k\in\mbox{dom}(s)\}\cap U\neq\emptyset$ as $\varphi(x)\in U$ and $x\supseteq s$ implies $\varphi(x)\in\bigcap\{x(n):n<\omega\}\subseteq\bigcap\{s(k):k\in\mbox{dom}(s)\}$. Let $x_V$ be an extension of $s$ such that $x_V(k)=V$ for $k\geqslant\mbox{dom}(s)$. Then $x_V\in\mathcal{U}$ by definition of $x_V$ and $\mathcal{U}$. Thus $x_V\in[s]\cap\mathcal{U}$ which contradicts the choice of $[s]$. As $U$ was arbitrary the Baire-measurability of $\varphi$ follows.

\end{enumerate}

Having established the above properties let us define $$\varphi^{-1}[\F]:=\left\{\varphi^{-1}[F]\subseteq X(\tau):F\in\F\right\}.$$ As a counterimage of a function preserves Boolean operations the family $\varphi^{-1}[\F]$ is a partition of the space $X(\tau)$ since $\F$ is a partition of $X$. Moreover $\varphi^{-1}[\F]$ consists of meager subsets in $X(\tau)$ by $2.$ above. We need only to show that for any subfamily $\F'\subseteq\F$ the set $\bigcup\{\varphi^{-1}[F]\subseteq X(\tau):F\in\F'\}$ has the Baire property in $X(\tau)$. For let $\F'\subseteq\F$ be arbitrary subfamily. As $\F$ is a Kuratowski partition of $X$ we have $\bigcup\F'=U\triangle M$ for some $U\in\tau$ and $M$ meager in $X$. Then
$$\bigcup\left\{\varphi^{-1}[F]\subseteq X(\tau):F\in\F'\right\}=\varphi^{-1}\left[\bigcup\F'\right]=\varphi^{-1}\left[U\triangle M\right]=\varphi^{-1}[U]\triangle\varphi^{-1}[M].$$

By $2.$ and $3.$ the last set has the Baire property in $X(\tau)$ \quad $\Box$

\medskip

\noindent{\bf Proof of Proposition 3:}

\medskip

Let $(X,\tau)$ be a Baire space with a Kuratowski partition. By {\bf Lemma 4} and {\bf Lemma 5} the space $X(\tau)$ is Baire metric space with a Kuratowski partition. \quad $\Box$


{\sc Ryszard Frankiewicz}\\ Institute of Mathematics, Polish
Academy of Sciences,\\ \'Sniadeckich 8, 00-950 Warszawa, Poland \\
{\em e-mail: rf@impan.pl}\\

{\sc S{\l}awomir Szczepaniak}\\ Institute of Mathematics, Polish
Academy of Sciences,\\ \'Sniadeckich 8, 00-950 Warszawa, Poland \\
{\em e-mail: szczepaniak@impan.pan.wroc.pl}

\end{document}